\documentclass{article}
\usepackage{amssymb}
\usepackage{amsfonts}
\usepackage{amsmath}

\input{tcilatex}
\begin{document}

\title{Characterization of the Pareto Social Choice Correspondence}
\author{Jerry S. Kelly \\
Department of Economics, Syracuse University\\
Syracuse, NY 13244-1020, USA\\
email: \texttt{jskelly@maxwell.syr.edu}}
\maketitle

\begin{abstract}
Independent, necessary and sufficient conditions are derived for a social
choice correspondence to be the one that selects exactly the Pareto optimal
alternatives.\medskip

\begin{itemize}
\item \noindent \textbf{Keywords:} Pareto $\cdot $ tops-in $\cdot $
balancedness $\cdot $ monotonicity $\cdot $ stability

\item \noindent \textbf{JEL:} D70 D71
\end{itemize}
\end{abstract}

\pagebreak

\begin{center}
{\LARGE Characterization of the Pareto Social Choice\smallskip }

{\LARGE Correspondence\footnote{%
The author is indebted to Shaofang Qi, Somdeb Lahiri, and an anonymous
reader for comments on earlier drafts.}}\bigskip
\end{center}

\section{Introduction\protect\medskip}

We characterize the social choice correspondence that, at each profile of
preferences, selects exactly the set of Pareto optimal alternatives. \ We
use one condition, balancedness, introduced in (Kelly and Qi 2019) and a
second, stability\ (related to a condition in Campbell, Kelly, and Qi 2018),
as well as tops-in and the Pareto condition. Although the collection of
Pareto optimal states at any situation has long been a central object of
concern in welfare economics, there is no work we know of characterizing
this mapping\footnote{%
Weymark (1984) characterizes the mapping to Pareto social \textit{preferences%
}. \ Here we are concerned with social choice correspondences that map to
social \textit{choice sets}.}. \ Our main theorem shows independent,
necessary, and sufficient conditions for the Pareto correspondence in the
case of five or more alternatives. \ Three other results, using weaker
conditions, for the cases of two, three, or four alternatives are also
included. \medskip

\section{Framework \protect\medskip}

Let $X$ with cardinality $|X|=m\geq 2$ be the finite set of \textbf{%
alternatives} and let $N=\{1,2,\ldots ,n\}$ with $n\geq 2$ be the set of 
\textbf{individuals}. A (strong) \textbf{ordering} on $X$ is a complete,
asymmetric, transitive relation on $X$ (so we exclude non-trivial individual
indifference). The top-ranked element of an ordering $r$ is denoted $r[1]$,
the next highest is denoted $r[2]$, etc. \ The set of all orderings on $X$
is $L(X)$. A \textbf{profile} $u$ is an element $(u(1),u(2),\ldots ,u(n))$
of the Cartesian product $L(X)^{N}$. \ If $x$ ranks above $y$ in $u(j)$, we
sometimes write\ $x\succ _{j}^{u}y$.\medskip

A \textbf{social choice correspondence} $G$ is a map from the domain\ $%
L(X)^{N}$ to non-empty subsets of $X$. \ One example that will play a role
here is the correspondence $T$ that maps profile $u$ to $T(u)$, the set of
all top-ranked alternatives at $u$:%
\begin{equation*}
T(u)=\{x\in X:\text{for some }i\in N\text{, }x=u(i)[1]\}
\end{equation*}

At profile $u$, alternative $x$ \textbf{Pareto dominates} $y$ if $x\succ
_{i}^{u}y$ for all $i$. \ Social choice correspondence $G$ satisfies the 
\textbf{Pareto condition} if for all $x$, $y$, and $u$, whenever $x$ Pareto
dominates $y$ at $u$, then $y\notin G(u)$. \ The \textbf{Pareto
correspondence}, $G_{P}$, is defined by%
\begin{equation*}
G_{P}(u)=\{x\in X:\text{there does not exist a }y\text{ that Pareto
dominates }x\text{ at }u\}
\end{equation*}%
\newline
Thus, the Pareto condition is

\begin{equation*}
G(u)\subseteq G_{P}(u)\text{ for all }u\text{.}
\end{equation*}

At profile $u$, alternative $x$ is \textbf{Pareto optimal} if there does not
exist an alternative $y$ that Pareto dominates $x$. \ Thus $G_{P}(u)$ is the
set of all Pareto optimal alternatives at $u$.\medskip

Clearly, the Pareto correspondence satisfies the Pareto condition, and this
is one of the conditions we use to characterize the Pareto correspondence. \
This may seem odd at first, but the condition of excluding dominated
alternatives is very weak and is used in characterizing many standard social
choice correspondences including the following which will appear several
times in this paper.\medskip

1. \ The \textbf{Borda correspondence}: \ Given profile $u$ and alternative $%
x$, the Borda score, $N(x,u)$, is the sum of the ranks in which $x$ appears
in the orderings in $u$. \ Alternative $x$ is a Borda winner at $u$ if $%
N(x,u)\leq N(y,u)$ for all $y$. \ The Borda correspondence selects at $u$
the set of all Borda winners at $u$. \medskip

2. \ The \textbf{plurality voting correspondence}: \ Given profile $u$ and
alternative $x$, the plurality score, $N(x,u)$, is the number of individuals
who have $x$ in the top rank. \ Alternative $x$ is a plurality winner at $u$
if $N(x,u)\geq N(y,u)$ for all $y$. \ The plurality correspondence selects
at $u$ the set of all plurality winners at $u$. \medskip

3. \ The \textbf{Copeland correspondence}: \ Given profile $u$ and
alternative $x$, the Copeland score, $N(x,u)$, is the number alternatives
defeated or tied by $x$ under simple majority voting. \ Alternative $x$ is a
Copeland winner at $u$ if $N(x,u)\geq N(y,u)$ for all $y$. \ The Copeland
correspondence selects at $u$ the set of all Copeland winners at $u$. \
\medskip \newline
The issue of using the Pareto condition is addressed again in the
conclusion.\medskip

The Pareto correspondence also satisfies other conditions commonly used in
social choice theory:\medskip

\textbf{Anonymity}: A social choice correspondence $G$ satisfies \textbf{%
anonymity} if, for every permutation $\rho $ on $N$, and every profile $u$%
,\medskip

\qquad \qquad $G(u(1),u(2),\ldots ,u(n))=G(u(\rho (1)),u(\rho (2),\ldots
,u(\rho (n))$.\medskip

\textbf{Neutrality}: Let $\theta $ be a permutation of $X$. \ If $%
S=\{x,y,\ldots ,z\}$ is a subset of $X$, we set $\theta (S)=\{\theta
(x),\theta (y),\ldots ,\theta (z)\}$. \ And if $R$ is an ordering on $X$, we
define $\theta (R)$ by $(\theta (x),\theta (y))\in \theta (R)$ if and only
if $(x,y)\in R$. Now we say a social choice correspondence $G$ satisfies 
\textbf{neutrality} if, for every permutation $\theta $ of $X$, and every
profile $u$,\medskip

\qquad \qquad $G(\theta (u(1)),\theta (u(2)),\ldots ,\theta (u(n)))=\theta
(G(u(1),u(2),\ldots ,u(n)))$. \medskip \newline
Those two additional conditions are not used in our characterization
theorems but are referred to in examples.\medskip

\section{Characterization with $m =2$\protect\medskip}

We start with a property used in all the characterizations in this
paper:\medskip

\qquad \qquad \textbf{Tops-in}: \ $T(u)\subseteq G(u)$ for all $u$ \medskip 
\newline
\relax Of course the tops rule choosing $T(u)$, the constant rule $G(u)=X$,
and the Pareto correspondence all satisfy tops-in. \ The plurality
correspondence, the Borda correspondence, and the Copeland correspondence do
not satisfy tops-in. \ Note that we are not restricting in this paper to
"desirable" conditions. \ Since the Pareto correspondence itself is not
desirable - if only because choice sets are often too large - any set of
characterizing conditions must include some that are undesirable.\medskip

\textbf{Theorem 1.} \ For $m=2$ and $n\geq 2$, let $G:L(X)^{N}\rightarrow
2^{X}\backslash \{\varnothing \}$ be a social choice correspondence
satisfying both of:\medskip

\qquad 1. \ The Pareto condition; and

\qquad 2. \ Tops-in;\medskip \newline
\relax then $G=G_{P}$, the Pareto correspondence.\medskip

\textbf{Proof:} \ The necessity of these conditions is obvious. \ Let $%
X=\{x,y\}$ and let $G$ be a social choice correspondence satisfying the
assumptions of the theorem.\medskip

\qquad \textbf{Case 1}. \ Only one of the two alternatives, say $x$, is in $%
T(u)$. \ Then by tops-in, $x\in G(u)$. \ By the Pareto condition, $y\notin
G(u)$. \ Therefore, $G(u)=\{x\}=G_{P}(u)$.\medskip

\qquad \textbf{Case 2}. \ Otherwise both alternatives, $x$ and $y$, are in $%
T (u)$ and then, by tops-in, $G (u) =\{x ,y\} =G_{P} (u)$. \ \ \ \ $\square $
\medskip

We now present examples showing that neither condition can be dropped in
Theorem 2:\medskip

\textbf{Example 1}. \ A rule other than $G_{P}$ satisfying tops-in (as well
as anonymity and neutrality) but not Pareto:\ $G(u)=X$ for all profiles $u$%
.\medskip

\textbf{Example 2}. \ A rule other than $G_{P}$ satisfying Pareto (as well
as anonymity and neutrality) but not tops-in:\ Let\ $m=2$ and $n\geq 3$ and
set\ $G(u)$ equal to the set of plurality winners at $u$.\medskip

\section{Characterization with $m =3$\protect\medskip}

When there are more than two alternatives, the properties of Theorem 1,
namely Pareto and tops-in, are \textit{not} sufficient to characterize the
Pareto correspondence.\medskip

\textbf{Example 3.} \ For $m\geq 3$, $G(u)=T(u)$ is distinct from $G_{P}(u)$
but satisfies Pareto and tops-in (as well as anonymity and
neutrality).\medskip

Note that $T(u)$ fails the following balancedness condition (Kelly and Qi
2019).\medskip

We say profile $v$ is \textbf{constructed from profile }$u$\textbf{\ by
transposition pair} $(\mathbf{x},\mathbf{y})$\textbf{\ via individuals }$i$%
\textbf{\ and }$j$ if at $u$, $x$ is immediately above $y$ for $i$ and $y$
is immediately above $x$ for $j$, and profile $v$ is just the same as $u$
except that alternatives $x$ and $y$ are transposed for $i$ and for $j$. A
social choice correspondence $G$ will be called \textbf{balanced} if, for
all $x$, $y$, $u$, $v$, $i$, and $j$, whenever profile $v$ is constructed
from $u$ by transposition pair $(x,y)$ via individuals $i$ and $j$, then $%
G(v)=G(u)$.\medskip \newline
The constant rule $G(u)=X$, the Pareto correspondence, Borda's rule, and the
Copeland correspondence all satisfy balancedness, but $T$ and the plurality
correspondence do not. In (Kelly and Qi, 2019) it is observed: \medskip

\begin{quotation}
[Balancedness] is a natural equity condition that simultaneously
incorporates some equal treatment for individuals short of anonymity, some
equal treatment of alternatives short of neutrality, and some equal
treatment of position of alternatives in orderings (for example, raising $x$
just above $y$ in the bottom two ranks for individual $j$ exactly offsets
lowering $x$ just below $y$ in the top two ranks for $i$).\medskip
\end{quotation}

We now show that, for three alternatives, incorporating balancedness with
tops-in and the Pareto condition forces $G=G_{P}$.\medskip

\textbf{Theorem 2}. \ For $m=3$ and $n\geq 2$, let $G:L(X)^{N}\rightarrow
2^{X}\backslash \{\varnothing \}$ be a social choice correspondence
satisfying all of:\medskip

\qquad 1. \ The Pareto condition;

\qquad 2. \ Tops-in;

\qquad 3. \ Balancedness; \medskip \newline
\relax then $G =G_{P}$, the Pareto correspondence. \medskip

\textbf{Proof}: The necessity of these conditions is obvious. \ \medskip

Now assume that $G$ satisfies all three conditions. \ We need to show that
if $w$ is Pareto optimal at $u$, then $w\in G(u)$. Suppose that Pareto
optimal $w\notin G(u)$. \medskip

Alternative $w$ can not be anyone's top at $u$ by tops-in. \ And $w$ can not
be everyone's bottom since it is Pareto optimal. \ So someone, say \#1, has $%
w$ in their second rank. \ Suppose that (1) $w$ is Pareto optimal, (2) $w$
is in \#1's second rank at $u$, (3) $w$ is in no one's top rank at $u$, and
(4) $w\notin G(u)$. In particular, assume that $x\neq z$ is at \#1's top.
Some individual, say \#2, has $w\succ _{2}^{u}x$ since $w$ is Pareto
optimal. For \#2, $w$ must be adjacent to $x$ with $x$ bottom-ranked.\medskip

\qquad \qquad 
\begin{tabular}{|l|l|l|}
\hline
1 & 2 & $\cdots $ \\ \hline
$x$ & $\vdots $ &  \\ 
$w$ & $w$ &  \\ 
$\vdots $ & $x$ & $\cdots $ \\ \hline
\end{tabular}
\medskip \newline
\relax By balancedness, $G(u^{\prime })=G(u)$ where $u^{\prime }$ is
obtained from $u$ by transposition pair $(x,w)$ for \#1 and \#2. \ So $%
w\notin G(u^{\prime })$. \ But that contradicts tops-in.\ \ \ \ $\square $ \
\medskip \medskip

To show the need for each condition in Theorem 2, we first observe that
Example 3 exhibits a rule other than $G_{P}$ satisfying all conditions of
Theorem 2 other than balancedness. \ Also, $G(u)=X$ for all profiles $u$
satisfies all conditions of Theorem 2 except Pareto. \medskip

\textbf{Example 4}. \ A rule other than $G_{P}$ satisfying all conditions of
Theorem 2 except tops-in: \ Let $m=n=3$ and let $u^{\ast }$ be a fixed
voter's paradox profile, say \medskip

\qquad \qquad 
\begin{tabular}{|l|l|l|}
\hline
1 & 2 & 3 \\ \hline
$x$ & $y$ & $z$ \\ 
$y$ & $z$ & $x$ \\ 
$z$ & $x$ & $y$ \\ \hline
\end{tabular}
\medskip \newline
\relax Define $G(u^{\ast })=\{x\}$ and $G(u)=G_{P}(u)$ for all $u\neq
u^{\ast }$. \ Balancedness is satisfied by $G$ because there are \textit{no}
transposition pairs at a voter's paradox profile. This rule is neither
anonymous nor neutral.\medskip

\section{Characterization with $m =4$\protect\medskip}

When there are more than three alternatives, the properties of Theorem 2,
namely Pareto, tops-in, and balancedness, are \textit{not} sufficient to
characterize the Pareto correspondence. \medskip

\textbf{Example 5}. \ Let $X=\;\{x,y,z,w\}$ with $n=3$. \ Consider fixed
profile $u^{\ast }$: \medskip

\qquad \qquad 
\begin{tabular}{|l|l|l|}
\hline
1 & 2 & 3 \\ \hline
$x$ & $y$ & $z$ \\ 
$y$ & $w$ & $w$ \\ 
$z$ & $x$ & $x$ \\ 
$w$ & $z$ & $y$ \\ \hline
\end{tabular}
\medskip \newline
\relax that has \textit{no} transposition pairs. Observe that $G_{P}(u^{\ast
})=\{x,y,z,w\}$. \medskip

Now define social choice correspondence $G$ as follows:\medskip

\qquad 1. \ $G (u^{ \ast }) =\{x ,y ,z\}$ $\;( =T (u^{ \ast }))$;

\qquad 2. \ For all other $u$, set $G(u)=G_{P}(u)$.\medskip \newline
\relax This correspondence (which fails anonymity and neutrality) satisfies
Pareto, tops-in, and balancedness, but differs from $G_{P}$ at $u^{\ast }$.
\ So we need to add some new condition to those of Theorem 2 in order to
characterize the Pareto correspondence for $m>3$. \ What won't work is using
anonymity and neutrality. \ Example 5 could be modified by constructing
subdomain $D$ consisting of all profiles obtained from the $u^{\ast }$ of
that example by permuting either $X$ or $N$ or both. \ Then define \medskip

\qquad 1. \ $G (u) =T (u)$ for all $u$ in $D$;

\qquad 2. \ $G(u)=G_{P}(u)$ otherwise.\medskip \newline
This $G$ satisfies Pareto, tops-in, balancedness, anonymity, and neutrality,
but differs from $G_{P}$ at every profile in $D$. \ So we introduce a new
property.\medskip

\textbf{Strong monotonicity}: \ A social choice correspondence $G$ satisfies 
\textbf{strong monotonicity\footnote{%
This is the same as correspondence monotonicity in Moulin (1983, p. 35).}}
if, for every $x\in X$, $i\in N$, and every profile $u$, if $x\in G(u)$, and
profile $u^{\prime }$ is constructed from $u$ by raising $x$ in $i$'s
ordering and leaving everything else unchanged, then

\begin{equation*}
x\in G(u^{\prime })\subseteq G(u)
\end{equation*}%
\newline
\relax Raising $x$ causes $x$ to be chosen again, but \textit{does not allow}
new alternatives to be chosen that weren't chosen before. The constant rule $%
G(u)=X$, Pareto, Borda, the plurality correspondence, and $T$ all satisfy
strong monotonicity. \medskip

Example 5, however, fails monotonicity in a significant way. \ If $z$, which
is in $G(u^{\ast })$, is raised one rank in profile $u^{\ast }$ for \#2,
then $G$ maps the resulting profile to $\{x,y,z,w\}$ and a new alternative, $%
w$, has been introduced to the choice set. \ We now show there is no way to
incorporate monotonicity without forcing $G=G_{P}$.\medskip

In this and the next section, we will want to show that when a social choice
correspondence $G$ satisfies certain properties it is $G_{P}$. \ If $G$
satisfies the Pareto condition, all alternatives that are \textit{not}
Pareto optimal at any profile $u$ are excluded from $G(u)$. \ What remains
is to show that every alternative that \textit{is} Pareto optimal at $u$ is
contained in $G(u)$. \ That leads us to consider the possibility of
alternatives $w$ in $G_{P}(u)\backslash G(u)$.\medskip

Given $G$, suppose that there exist profiles where some Pareto optimals are
not chosen by $G$. Consider the non-empty collection $\mathbb{C}\subseteq
L(X)^{N}$ of all profiles for which there exists at least one Pareto optimal
alternative that is not chosen.\ \ For profile $v\in \mathbb{C}$, when there
exists at least one individual $i$ and alternative $w$ such that $w=v(i)[t]$
and $w\in G_{P}(v)\backslash G(v)$ but there does not exist an $s<t$, an
individual $j$ and alternative $y$ such that $y=v(j)[s]$ and $y\in
G_{P}(v)\backslash G(v)$, set the \textbf{height} at $v$ as $h(v)=t$. \
While $h(u)$ is uniquely determined from $G$ and $u$, that is not true for
the specification of individuals and alternatives. \ For example, a profile $%
u$ might have $h(u)=2$, both because \#1 has unchosen but Pareto optimal $x$
in second rank and \#2 has unchosen but Pareto optimal $y$ in second rank.
Of all these values of $h(v)$ for profiles $v$ in $\mathbb{C}$, let $H(G)$,
the \textbf{height of }$G$, be the \textit{minimum} $h(v)$ (corresponding to
a \textit{highest} ranked alternative from $G_{P}(v)\backslash G(v)$%
).\medskip 

Suppose that for correspondence $G$, the height $H(G)$ is defined. \
Consider the possibility that $H(G)=m$. \ This requires in part that there
is a profile $u$ and alternative $x$ such that:

\qquad (1) $x$ is Pareto optimal at $u$;

\qquad (2) $x$ is not chosen;

\qquad (3) everyone has $x$ in rank $m$.

But (3) contradicts (1). \ Therefore,

\begin{equation*}
H(G)\leq m-1\text{.}
\end{equation*}%
\newline

With $n>2$, the plurality correspondence has height 1 as does the Copeland
correspondence. For $m>2$, height $H(G)>1$ if and only if $G$ satisfies
tops-in\footnote{%
If $m=2$, then $H(G)>1$ is not possible.}. The correspondence $T$ has height
2. \ Height is not defined for $G_{P}$ or the constant rule $G(u)=X$, since
those correspondences always choose all Pareto optimals.\medskip 

Height is used in this section just as an organizational device in the proof
of Theorem 3. \ But in the next section it is central to the proof of
Theorem 4. \medskip

\textbf{Theorem 3}. \ For $m =4$ and $n \geq 2$, let $G :L (X)^{N}
\rightarrow 2^{X}\backslash \{\varnothing \}$ be a social choice
correspondence satisfying all of: \medskip

\qquad 1. \ The Pareto condition;

\qquad 2. \ Tops-in;

\qquad 3. \ Balancedness;

\qquad 4. \ Strong monotonicity; \medskip \newline
\relax then $G=G_{P}$, the Pareto correspondence. \medskip

\textbf{Proof}: The necessity of these conditions is obvious. \medskip

Now assume that $G$ satisfies all four conditions. We need to show that if $%
w $ is Pareto optimal at $u$, then $w\in G(u)$. \ Assume otherwise; then $%
H(G)$ is defined. \ Height $H(G)=1$ is excluded by tops-in and $H(G)=4$ is
excluded by the Pareto condition. \ Two possibilities then remain. \medskip

\textbf{Case 1}. \ $H(G)=2$. \ Suppose that $w\in G_{P}(u)\backslash G(u)$
is in \#1's second rank at $u$ given in part by: \medskip

\qquad \qquad 
\begin{tabular}{|l|l|l|}
\hline
1 & 2 & $\cdots $ \\ \hline
$x$ & $\vdots $ &  \\ 
$w$ &  &  \\ 
$\vdots $ & $w$ & $\cdots $ \\ 
& $\vdots $ &  \\ 
& $x$ &  \\ \hline
\end{tabular}
\medskip \newline
\relax where $x$ is at \#1's top (so $x\in G(u)$ by tops-in). \ Some
individual, say \#2, has $w\succ _{2}^{u}x$ since $w$ is Pareto optimal. \
Construct $u^{\prime }$ from $u$ by raising $x$ up to just below $w$ for
\#2; $w$ remains Pareto optimal. Now $x\in G(u^{\prime })\subseteq G(u)$ by
monotonicity. \ In particular, $w\notin G(u^{\prime })$. \ By balancedness, $%
G(u^{\prime \prime })=G(u^{\prime })$ where $u^{\prime \prime }$ is obtained
from $u^{\prime }$ by transposition pair $(x,w)$ via \#1 and \#2. \ So $%
w\notin G(u^{\prime \prime })$. \ But that contradicts tops-in since, at $%
u^{\prime \prime }$, \#1 now has $w$ top-ranked.\medskip

\textbf{Case 2}. \ $H(G)=3$. \ Since $w$ is Pareto optimal, it must be
ranked higher than each of the other alternatives; without loss of
generality, $u$ is, in part: \medskip

\qquad \qquad 
\begin{tabular}{|l|l|l|l|}
\hline
1 & 2 & 3 & $\cdots $ \\ \hline
&  &  &  \\ 
$\vdots $ & $\vdots $ & $\vdots $ &  \\ 
&  &  & $\cdots $ \\ 
$w$ & $w$ & $w$ &  \\ 
$x$ & $y$ & $z$ &  \\ \hline
\end{tabular}
\medskip \newline
\relax Take the alternative just above $w$ for \#1 to be say, $y$. \ Then
construct $u^{\ast }$ from $u$ by transposition pair $(y,w)$ via \# 1 and
\#2. \ By balancedness, $w\notin G(u^{\ast })$ where $w$ is Pareto optimal
at $u^{\ast }$ and ranked second by \#1 contrary to our assumption that $%
H(G)=3$. \ \ \ \ $\square $ \medskip

For examples showing the need for each condition in Theorem 3 (with $n=3$
and $m=4$), we first observe that Example 5 exhibits a rule different from $%
G_{P}$ satisfying all conditions of Theorem 3 other than monotonicity. And $%
G(u)=X$ for all profiles $u$ is a rule different from $G_{P}$ satisfying all
conditions of Theorem 3 except Pareto. (Another: set $G(u)=G_{P}(u)$ except
at profiles $u$ where everyone has the same top and the same second-ranked
alternative, at such profiles, $G(u)$ is the set consisting of those two
alternatives.) \ For the others: \medskip

\textbf{Example 6}. \ A rule different from $G_{P}$ satisfying all
conditions of Theorem 3 except tops-in: \ Fix one alternative $t$, and then
set $G(u)=G_{P}(u)\backslash \{t\}$ unless $G_{P}(u)=\{t\}$, in which case
set $G(u)=G_{P}(u)=\{t\}$. \medskip

\textbf{Example 7}. \ A rule different from $G_{P}$ satisfying all
conditions of Theorem 3 except balancedness: \ $G(u)=T(u)$. \medskip

A variant of Theorem 3 will appear as Theorem 5 near the end of the next
section.\medskip

\section{Characterization with $m \geq 5$\protect\medskip}

When there are more than four alternatives, the properties of Theorem 3,
namely Pareto, tops-in, balancedness, and strong monotonicity, are \textit{%
not} sufficient to characterize the Pareto correspondence. \medskip

\textbf{Example 8}. \ Let $X=\{x,y,z,w,t\}$ and $n=2$. \ Define social
choice correspondence $G$ as follows. \ First we identify a subdomain $D$ of 
$L(X)^{N}$ that consists of just the two fixed profiles \medskip

\qquad \qquad $u$:\qquad 
\begin{tabular}{|l|l|}
\hline
1 & 2 \\ \hline
$x$ & $z$ \\ 
$y$ & $t$ \\ 
$w$ & $w$ \\ 
$z$ & $x$ \\ 
$t$ & $y$ \\ \hline
\end{tabular}%
\qquad and $u^{\ast }$:\qquad 
\begin{tabular}{|l|l|}
\hline
1 & 2 \\ \hline
$z$ & $x$ \\ 
$t$ & $y$ \\ 
$w$ & $w$ \\ 
$x$ & $z$ \\ 
$y$ & $t$ \\ \hline
\end{tabular}
\medskip \newline
\relax At these profiles in $D$, set $G(u)=G(u^{\ast })=\{x,z\}$, the
top-most alternatives (thus $w$ is not chosen even though it is Pareto
optimal at these profiles). \ For all profiles $v$ in $L(X)^{N}\backslash D$%
, set $G(v)=G_{P}(v)$. \medskip

Clearly $G$ satisfies the Pareto condition and tops-in. \ For balancedness,
observe that there does not exist a transposition pair at either of the
profiles in $D$. \ Accordingly, if $v$ is obtained from $u$ by pairwise
transposition, both $u$ and $v$ are in $L(X)^{N}\backslash D$ where $%
G(u)=G_{P}(u)$ and $G(v)=G_{P}(v)$. \ Since $G_{P}$ satisfies balancedness,
so does $G$. \medskip

All that remains is monotonicity. \ If $v$ and $u$ are both in $%
L(X)^{N}\backslash D$, and $v$ arises from $u$ by raising a chosen
alternative $x$, then, because $G=G_{P}$ there, and $G_{P}$ satisfies
monotonicity, we cannot have a violation of monotonicity by $G$. \ Neither
profile in $D$ can arise from the other by raising a chosen alternative. \
So, if $G$ fails monotonicity, it has to be because raising a chosen
alternative takes you from $D$ into $L(X)^{N}\backslash D$\ or from $%
L(X)^{N}\backslash D$ into $D$ . \medskip

From $D$ to $L(X)^{N}\backslash D$: \ Suppose that $v$ is constructed by
raising $x$ at $u$ where $G(u)=\{x,z\}$ (all other cases are dealt with by
simple analogs of this argument). \ This must be for individual \#2 and
raising $x$ means that $x$ will now Pareto dominate $w$. \ Since $y$ and $t$
remain Pareto dominated, $G(v)=G_{P}(v)$ will be $\{x,z\}$ or $\{x\}$ (if $x$
is raised to \#2's top). \ In either case, we have \medskip

\qquad \qquad \qquad \qquad $x \in G (v) \subseteq G (u)$. \medskip

From $L(X)^{N}\backslash D$ to $D$: Suppose that $u\in D$ is constructed by
raising an alternative from a profile $q$ in $L(X)^{N}\backslash D$. \ So
construct $q$ from $u$ by lowering an alternative for someone. \ If $y$ or $%
t $ is lowered, it remains Pareto dominated and so not in $G(q)$. \ If $w$
is lowered, it becomes Pareto dominated and so not in $G(q)$. \ If $x$ or $z$
is lowered and is chosen at $G(q)$, then we have e.g., \medskip

\qquad \qquad \qquad \qquad $x\in G(u)\subseteq G(q)$. \medskip \newline
\relax In each case, strong monotonicity is confirmed. Note that $G$, while
not neutral, is anonymous ($u^{\ast }$ is obtained from $u$ by a permutation
on $N=\{1,2\}$). \medskip

Our next example modifies Example 8 by allowing $n>2$. \medskip

\textbf{Example 9}. Again we identify a subdomain $D$ of $L(X)^{N}$. \ We
start from a list $C$ of eight orderings on $X=\{x,y,z,w,t\}$: \medskip

\qquad 
\begin{tabular}{|l|}
\hline
$x$ \\ 
$y$ \\ 
$w$ \\ 
$z$ \\ 
$t$ \\ \hline
\end{tabular}%
\qquad 
\begin{tabular}{|l|}
\hline
$z$ \\ 
$t$ \\ 
$w$ \\ 
$x$ \\ 
$y$ \\ \hline
\end{tabular}%
\qquad 
\begin{tabular}{|l|}
\hline
$x$ \\ 
$y$ \\ 
$z$ \\ 
$t$ \\ 
$w$ \\ \hline
\end{tabular}%
\qquad 
\begin{tabular}{|l|}
\hline
$x$ \\ 
$z$ \\ 
$y$ \\ 
$t$ \\ 
$w$ \\ \hline
\end{tabular}%
\qquad 
\begin{tabular}{|l|}
\hline
$x$ \\ 
$z$ \\ 
$t$ \\ 
$y$ \\ 
$w$ \\ \hline
\end{tabular}%
\qquad 
\begin{tabular}{|l|}
\hline
$z$ \\ 
$x$ \\ 
$y$ \\ 
$t$ \\ 
$w$ \\ \hline
\end{tabular}%
\qquad 
\begin{tabular}{|l|}
\hline
$z$ \\ 
$x$ \\ 
$t$ \\ 
$y$ \\ 
$w$ \\ \hline
\end{tabular}%
\qquad 
\begin{tabular}{|l|}
\hline
$z$ \\ 
$t$ \\ 
$x$ \\ 
$y$ \\ 
$w$ \\ \hline
\end{tabular}
\medskip \newline
\relax The first two orderings in $C$ are the orderings in the profiles in
the subdomain $D$ of Example 8. \ The next six orderings consist of the six
possible ways of ordering $\{x,y,z,t\}$ subject to $x\succ y$ and $z\succ t$
with $w$ then appended at the bottom. \ Subdomain $D$ consists of all
profiles made up of just these orderings subject to the condition that 
\textit{each of the first two orderings occurs exactly once}. \ For any
profile $u\in D$, we set $G(u)=T(u)=\{x,z\}$. For all profiles $v$ in $%
L(X)^{N}\backslash D$, set $G(v)=G_{P}(v)$.\medskip

Clearly $G$ satisfies the Pareto condition and tops-in. \ It is
straightforward to check that balancedness and strong monotonicity hold. \
Note that $G$ is anonymous (as any permutation of $N$ takes a profile in $D$
to another profile in $D$).\medskip

We next illustrate the need for the restriction that each of the first two
orderings occurs exactly once. \ First, we need both to occur so that $w$ is
Pareto optimal. \ To illustrate the need for the orderings to not occur more
than once, suppose that $n=3$ and consider the domain $D^{\ast }$ of
profiles made up of the eight orderings above (\textit{without} the
restriction that each of the first two orderings occurs exactly once). \
Define the correspondence $G^{\ast }$ by $G^{\ast }(u^{\ast })=T(u^{\ast })$
if $u^{\ast }\in D^{\ast }$ and $G^{\ast }(u)=G_{P}(u)$ if $u\notin D^{\ast
} $. \ Then look at the profiles \medskip

\qquad \qquad $u$: \ \ 
\begin{tabular}{|l|l|l|}
\hline
1 & 2 & 3 \\ \hline
$x$ & $z$ & $z$ \\ 
$y$ & $t$ & $t$ \\ 
$w$ & $w$ & $w$ \\ 
$z$ & $x$ & $x$ \\ 
$t$ & $y$ & $y$ \\ \hline
\end{tabular}%
\ \ and $v$: \ \ 
\begin{tabular}{|l|l|l|}
\hline
1 & 2 & 3 \\ \hline
$x$ & $z$ & $z$ \\ 
$y$ & $t$ & $t$ \\ 
$w$ & $x$ & $w$ \\ 
$z$ & $w$ & $x$ \\ 
$t$ & $y$ & $y$ \\ \hline
\end{tabular}
\medskip \newline
\relax$G(u)=\{x,z\}$ since $u\in D^{\ast }$; then $v$ is obtained by raising 
$x$ above $w$ for \#2. \ But $G(v)=G_{P}(v)=\{x,z,w\}$, a violation of
strong monotonicity. \medskip

Thus, to achieve a characterization of the Pareto correspondence for $m\geq
5 $, we need another new condition. \ We introduce a correspondence analog
of a property of social choice \textit{functions }that are resolute, i.e.,
single alternatives are chosen at every profile (Campbell, Kelly, and Qi,
2018). \ A social choice \textit{function} $g$ satisfies \textbf{stability}
if for every pair of profiles, $u$ and $u^{\ast }$, and for every pair $x$, $%
y$ of alternatives, and every individual $i$, if $x=g(u)$ and, in ordering $%
u(i)$, alternative $y$ is adjacent to and just below $x$, then if $u^{\ast }$
is obtained from $u$ by \textit{only} lowering $x$ to just below $y$ for $i$%
, we must have either $g(u^{\ast })=x$ or $g(u^{\ast })=y$. \ A small change
in the profile results in a restricted set of possible outcomes. \ Here,
where we are dealing with correspondences, we consider the same alteration
of profiles, but allow somewhat different consequences, but still with a
focus on "small changes cause small effects."\footnote{%
In (Campbell, Kelly, and Qi, 2018), it is observed that in Richard Bellman's
autobiography [1984, p. 181], he writes: "Change one small feature, and the
structure of the solution was strongly altered. There was no stability!"}. \
A social choice \textit{correspondence} $G$ satisfies \textbf{strong
stability} if for every pair of profiles, $u$ and $u^{\ast }$, and for every
pair $x$, $y$ of alternatives, and every individual $i$, if $x\in G(u)$ and,
in ordering $u(i)$, alternative $y$ is adjacent to and just below $x$, then
if $u^{\ast }$ is obtained from $u$ by \textit{only} lowering $x$ to just
below $y$ for $i$, we must have \textit{exactly one} of the following
outcomes hold at $u^{\ast }$: \medskip

\qquad \qquad a. \ $G (u^{ \ast }) =G (u)$;

\qquad \qquad b. \ $G(u^{\ast })=G(u)\backslash \{x\}$ where $x\in G(u)$;

\qquad \qquad c. \ $G (u^{ \ast }) =G (u) \cup \{y\}$ where $y \notin G (u)$%
. \medskip \newline
\relax If $G (u^{ \ast })$ differs from $G (u)$, it must either drop $x$ or
add $y$, \textit{but not both }(stability here is "strong" because of the
"not both" requirement). \medskip

The Pareto correspondence $G_{P}$ satisfies the strong stability property. \
Consider profile $u$ with $x\in G_{P}(u)$, where $y$ is just below $x$ in $%
u(i)$, and $u^{\ast }$ differs from $u$ only in that $x$ is brought down
just below $y$ in $i$'s ordering. \ If $z$ is any element of $X\backslash
\{x,y\}$, then $z$ is dominated by an element at $u^{\ast }$ if and only if
it is dominated by that same element at $u$. \ Hence $z\in G_{P}(u^{\ast })$
if and only if $z\in G_{P}(u)$. \ So $G_{P}(u^{\ast })$ can differ from $%
G_{P}(u)$ only by losing $x$ or gaining $y$. \ But if $y$ is gained, then $%
y\notin G_{P}(u)$. \ It must be that $y$ was Pareto dominated by $x$ at $u$.
\ But then at $u^{\ast }$, all individuals other than $i$ still prefer $x$
to $y$ and so $x$ is still Pareto optimal at $u^{\ast }$, i.e., $G_{P}(u)$ 
\textit{can not} both lose $x$ and gain $y$. \medskip

Clearly $G$ given by $G (u) =X$ for all $u$ satisfies strong stability. \ So
does the following rule: \ Set $G (u) =X$ unless, at $u$, there is a common
alternative $t$ at everyone's bottom rank; then set $G (u) =X\backslash
\{t\} $. \medskip

Many rules fail strong stability. Dictatorship fails strong stability. \ $%
G(u)=T(u)$, the Borda correspondence, and the plurality correspondence all
fail strong stability. \ Also, the correspondences of Examples 8 and 9 
\textit{do not} satisfy strong stability. For example, in Example 8,
consider profile $u\in D$: \medskip

\qquad \qquad 
\begin{tabular}{|l|l|}
\hline
1 & 2 \\ \hline
$x$ & $z$ \\ 
$y$ & $t$ \\ 
$w$ & $w$ \\ 
$z$ & $x$ \\ 
$t$ & $y$ \\ \hline
\end{tabular}
\medskip \newline
\relax where $G (u) =\{x ,z\}$. \ Construct $u^{ \ast } \notin D$ from $u$
by bringing $x$ down just below $y$ for individual 1. $G (u^{ \ast }) =G_{P}
(u) =\{x ,y ,z ,w\}$, a violation of strong stability. \medskip \ \ 

That strong stability fails for plurality and $T$ is seen at profile $u$:
\medskip

\qquad \qquad 
\begin{tabular}{|l|l|l|}
\hline
1 & 2 & 3 \\ \hline
$x$ & $y$ & $z$ \\ 
$a$ & $\vdots $ & $\vdots $ \\ 
$w$ &  &  \\ 
$\vdots $ &  &  \\ \hline
\end{tabular}
\medskip \newline
\relax Bringing $x$ down just below a for \#1 both loses $x$ and gains $a$.
\ Borda and Copeland also fail strong stability; consider profile $u$:
\medskip

\qquad \qquad 
\begin{tabular}{|l|l|l|}
\hline
1 & 2 & 3 \\ \hline
$x$ & $x$ & $y$ \\ 
$y$ & $y$ & $x$ \\ 
$\vdots $ & $\vdots $ & $\vdots $ \\ \hline
\end{tabular}
\medskip \newline
\relax Bringing $x$ down just below $y$ for \#1 both loses $x$ and gains $y$%
. \medskip

Note that strong stability does \textit{not} imply strong monotonicity.
Consider the correspondence that always selects individual \#1's
bottom-ranked alternative. Strong monotonicity treats raising a \textit{%
chosen} $y$ above anything. Strong stability only treats raising (any) $y$
above an adjacent \textit{chosen} $x$.\medskip

We will employ strong stability to complete a characterization of the Pareto
correspondence with five or more alternatives.\medskip

Suppose $G$ is a correspondence satisfying balancedness. \ Let $x$ be an
alternative, $u$ a profile with $x$ Pareto optimal but $x\notin G(u)$, and $%
i $ an individual with $x$ in rank $k$ where $x\notin G(u)$. The following
two propositions are inconsistent:\medskip

\qquad A. \ $H(G)=k$;\medskip

\qquad B. \ There is an alternative $\alpha $ adjacent to and just above $x$
for individual $i$ and adjacent to and just below $x$ for some other
individual $j$.\medskip \newline
If both were true at $u$, construct profile $u^{\ast }$ by transposing $x$
and $\alpha $ for individuals $i$ and $j$. \ Then $x$ is also Pareto optimal
at $u^{\ast }$ and, by balancedness, $x\notin G(u^{\ast })$. \ But then $%
H(G)<k$.\medskip

We will use this fact in two different ways:\medskip

\qquad (i) \ If we have assumed $H(G)=k$, then we must exclude the existence
of a profile with an alternative $\alpha $ adjacent to and just above $x$
(with $x$ in rank $k$) for individual $i$ and adjacent to and just below $x$
for some other individual $j$.\medskip

\qquad (ii) If we have a profile for which there exists an alternative $%
\alpha $ adjacent to and just above Pareto optimal but not chosen $x$ (in
rank $k$) for individual $i$ and adjacent to and just below $x$ for some
other individual $j$, we must reject $H(G)=k$.\medskip \newline
Either of these two uses will be called a \textbf{height argument}.\medskip

Note that in the following we do not assume monotonicity.\medskip

\textbf{Theorem 4}. \ For $m\geq 5$ and $n\geq 2$, let $G:L(X)^{N}%
\rightarrow 2^{X}\backslash \{\varnothing \}$ be a social choice
correspondence satisfying all of: \medskip

\qquad 1. \ The Pareto condition;

\qquad 2. \ Tops-in;

\qquad 3. \ Balancedness;

\qquad 4. \ Strong stability;\medskip \newline
\relax then $G=G_{P}$, the Pareto correspondence. \medskip

\textbf{Proof}: The necessity of these conditions is obvious. \medskip

Now assume that $G$ satisfies all four conditions but $G\neq G_{P}$. \ Then $%
H(G)$ is defined; say $H(g)=k$ for $1\leq k\leq m$. \ So there exists a
profile $u$ and alternative $w$ such that $w$ is Pareto optimal, $w$ is not
chosen, and $w$ is in the $k$th rank for, say, individual \#1. \ We will
show that each possible value of $k$ leads to a contradiction. \ Some values
of $H(G)$ are easily dealt with. \ $H(G)=1$ violates tops-in while $H(G)=m$
violates the Pareto condition. \ But we seek a more general, systematic
approach to deal will all possible values of $k$.\medskip

Without loss of generality, suppose profile $u$ is given in part by\medskip

\qquad \qquad $u$: 
\begin{tabular}{|l|l|l|}
\hline
1 & 2 & $\cdots $ \\ \hline
$a$ &  &  \\ 
$b$ &  &  \\ 
$\vdots $ &  &  \\ 
$c$ &  & $\cdots $ \\ 
$w$ &  &  \\ 
$d$ &  &  \\ 
$\vdots $ &  &  \\ 
$e$ &  &  \\ \hline
\end{tabular}%
\medskip \newline
with $w$ in the $k$th rank for \#1. \ We will methodically work down the top
set of alternatives in \#1's ranking in order from $a$ to $c$.\medskip

We know $a\in G(u)$ by tops-in.\medskip

\textbf{Case 1}. \ Some $i>1$ has $a\succ _{i}^{u}w$. \ \medskip

Without loss of generality, we take this to be \#2\medskip

\qquad \qquad $u$: 
\begin{tabular}{|l|l|l|}
\hline
1 & 2 & $\cdots $ \\ \hline
$a$ &  &  \\ 
$b$ &  &  \\ 
$\vdots $ & $p$ &  \\ 
$c$ & $a$ & $\cdots $ \\ 
$w$ & $\vdots $ &  \\ 
$d$ & $w$ &  \\ 
$\vdots $ &  &  \\ 
$e$ &  &  \\ \hline
\end{tabular}%
\medskip \newline
\textbf{Step 1}. \ If $a$ is at \#2's top, go to Step 2. Otherwise there is
a $p$ just above $a$ for \#2; construct $u^{\prime }$ from $u$ by raising $a$
just above $p$. \ By tops-in, $a$ is still chosen. Also $w$ is still Pareto
optimal and not chosen (otherwise going from $u^{\prime }$ to $u$ would
violate strong stability). \ Do this until $a$ is at the top of \#2's order
with $w$ still Pareto optimal but not chosen.\medskip \newline
\textbf{Step 2}. \ If, for \#1, $a$ is not adjacent to $w$, bring $a$ down
one rank. \ Then $a$ is still chosen by \#2 having $a$ at the top with $w$
still Pareto optimal but not chosen (otherwise going back would violate
strong stability). \ Do this until $a$ is adjacent and just above $w$ for
\#1. Call the resulting profile $u^{\ast }$.\medskip \newline
\textbf{Step 3}. \ Since $w$ is Pareto optimal, someone, say \#3, must have $%
w\succ _{3}^{u^{\ast }}a$. \ If $a$ is not adjacent to $w$ for \#3, raise $a$
one rank to get profile $u^{\ast \prime }$. \ By tops-in, $a$ is still
chosen. \ Continue in this way until $a$ is adjacent to and just below $w$
for \#3 with $w$ still Pareto optimal but not chosen (otherwise going back
would violate strong stability). \ Then a height argument for $a$ and $w$
for \#1 and \#3 shows $H(G)=k$ is violated.\medskip

\textbf{Case 2}. \ For all $i>1$, $w\succ _{i}^{u}a$. \ \medskip

For this case, we shift attention to $b$, the next lower alternative in
\#1's ranking. \ If $k=2$, go to the argument right after Subcase 2-2B. 
Assuming then $k\geq 3$, so $b$ is above $w$, there are several
possibilities:\medskip 

\textbf{Subcase 2-1}. \ $b$ is not Pareto optimal. \ Only $a$ can Pareto
dominate $b$ and, since for all $i>1$, $w\succ _{i}^{u}a$, we also have all $%
i>1$, $w\succ _{i}^{u}b$. \medskip

\textbf{Subcase 2-2}. \ If $b$ is Pareto optimal then because \#1 ranks $b$
above position $H(G)$, alternative $b$ must be chosen. \ Then we divide our
analysis just as we did for $a$.\medskip

\textbf{Subcase 2-2A}. \ Suppose some $i>1$ has $b\succ _{i}^{u}w$. \
Without loss of generality, we take this to be \#2:\medskip

\qquad \qquad $u$: 
\begin{tabular}{|l|l|l|}
\hline
1 & 2 & $\cdots $ \\ \hline
$a$ &  &  \\ 
$b$ &  &  \\ 
$\vdots $ &  &  \\ 
$c$ & $b$ & $\cdots $ \\ 
$w$ & $\vdots $ & $w$ \\ 
$d$ & $w$ & $\vdots $ \\ 
$\vdots $ & $\vdots $ & $a$ \\ 
$e$ & $a$ &  \\ \hline
\end{tabular}%
\medskip \newline
If $b$ is not \#2's top, raise $b$ one rank. \ Then $b$ is still Pareto
optimal and $b$ is chosen because \#1 ranks $b$ above $H(G)$. \ Alternative $%
w$ is still Pareto optimal. \ More, $w$ is not chosen; for if it were then
bringing chosen $b$ back down one rank would violate strong stability. \
Continue in this way, raising $b$ to \#2's top.\medskip

Then, just as in Step 2 above, bring $b$ down for \#1 until it is just above 
$w$. \ At the resulting profile, $b$ is chosen by tops-in (via \#2) while $w$
is Pareto optimal but not chosen. \ Since $w$ is Pareto optimal, someone,
say \#3, has $b$ below $w$. \ If $b$ is not adjacent to $w$ for \#3, raise $%
b $ one rank to get profile $u^{\ast \prime }$ (by successive applications
of strong stability). \ Since $b$ is Pareto optimal and \#1 ranks $b$ above
position $H(G)$, we see $b$ is still chosen. \ Continue in this way until $b$
is adjacent to and just below $w$ for \#3. \ Then a height argument shows $%
H(G)=k$ is violated.\medskip

\textbf{Subcase 2-2B}. \ Otherwise for all $i>1$, $w\succ _{i}^{u}b$.\medskip

Continue in this fashion, for each alternative $t$ that \#1 ranks above $w$.
\ For each alternative we find a contradiction from a height argument or we
learn that for everyone other than \#1 that alternative is ranked below $w$.
\ This procedure stops in only a very few possible ways.\medskip 

\qquad (1) If some alternative $t$ ranked above $w$ by \#1 is chosen and
someone else has $t$ ranked above $w$, then we move $t$ around in ways that
lead to a height argument contradicting $H(G)=k$. \ In particular this has
to happen if $H(G)-1>m-H(G)$.\medskip 

\qquad (2) If $H(G)-1=m-H(G)$ and we haven't run into a contradiction from a
height argument, then every $i>1$ has all of $\{a,b,...,c\}$ ranked below $w$
and everything in $\{d,...,e\}$ ranked above $w$ with everyone ranking $w$
in position $k$. \ Then \#2's top, say $d$, will be ranked above $w$ by
someone else (when $n\geq 3$; the case $n=2$ can be dealt with
straightforwardly) and then an argument like that for $a$ will lead to a
contradiction of $H(G)=k$.\medskip

\qquad (3) If $H(G)-1<m-H(G)$ and no inconsistency had been found,\ raise $a$
up to just below $w$ for some $i>1$. \ Next lower $a$ to just above $w$ for
\#1 as before; at each step in this process, alternative $a$ is ranked above
all of $\{b,...,c\}$ by individual $i$ and ranked above all of $\{w,d,...,e\}
$ for \#1, so $a$ always remains Pareto optimal. Since, for \#1, $a$ appears
above the $k$th rank, it can not be chosen or $H(G)=k$ would be violated.
Finally apply a height argument. \ \ \ \ $\square $\medskip 

Regarding the need for each condition, Example 9 shows a rule different from 
$G_{P}$ satisfying all conditions of Theorem 4 except strong stability. A
rule different from $G_{P}$ satisfying all conditions of Theorem 4 except
Pareto is:\ $G(u)=X$ for all profiles $u$. For tops-in consider the next
example.\medskip

\textbf{Example 10}. \ $m\geq 5$, $n>m$. \ $G(u)=G_{P}(u)\backslash \{x:$
for some $y$, exactly $n-1$ of the $n$ individuals prefer $y$ to $x\}$. \
(Since $n>m$, at profile $u$ there is at least one alternative at the top
for at least two individuals. \ This alternative is Pareto optimal and can't
lose $(n-1)$-to-$1$ to any alternative. \ Thus $G(u)$ is non-empty.) This
clearly satisfies Pareto (and anonymity and neutrality). \ It fails tops-in
as $x$ could be \#1's top and everyone else's bottom and thus lose in some $%
(n-1)$-to-$1$ votes. \ For balancedness, construction of $u^{\ast }$ by an
interchange of transposition pairs from $u$ isn't possible if one Pareto
dominates the other and doesn't change any $(n-1)$-to-$1$ vote; $G(u^{\ast
})=G(u)$. \ For stability, suppose that $u^{\ast }$ is constructed from $u$
by bringing $x$ down just below adjacent $y$. \ Can we have $G(u^{\ast
})=G(u)\backslash \{x\}$ where $x\in G(u)$?\ \ It can't be because now $y$
Pareto dominates $x$ for then previously $y$ beat $x$ by $(n-1)$-to-$1$ and $%
x\in G(u)$ would not be true. \ But it could be that at $u$, alternative $x$
lost to $y$ by $(n-1)$-to-$1$ but now only $n-2$-to-$2$. But then for $y$ to
move into $G(u)$ it would have to have been excluded at $u$. This requires a 
$z\neq x$ such that at least $n-1$ individuals prefer $z$ to $y$ at $u$; but
then this must also be true at $u^{\ast }$ so $y$ can't move into $G(u^{\ast
})$.\medskip

For balancedness, consider the next correspondence which is also anonymous
and neutral:\medskip

\textbf{Example 11}. \ Assume that $n=4$ and $m=5$ (or more generally, $%
n=m-1 $).\medskip

Partition $L(X)^{N}$ as $D\cup D^{\ast }$ where $D^{\ast }$ is the set of
profiles satisfying the following three properties:\medskip

\qquad 1. \ $u\in D^{\ast }$ implies there is a (unique) alternative $\phi
(u)$ that in

\qquad \qquad everyone's ordering at $u$ has $\phi (u)$ in the next-to-last
rank;

\qquad 2. \ For every alternative in $w\backslash \{\phi (u)\}$, one person
has that

\qquad \qquad alternative ranked at the bottom (just below $\phi (u)$);

\qquad 3. \ Every alternative is Pareto optimal at $u$ (note $\phi (u)$ is
already

\qquad \qquad guaranteed to be Pareto optimal by property 2);\medskip 
\newline
and $D=L(X)^{N}\backslash D^{\ast }$.\medskip

Define collective choice correspondence $G$ as follows:\medskip

\qquad If $u\in D$, then $G(u)=G_{P}(u)$, the Pareto optimals at $u$;

\qquad If $u\in D^{\ast }$, then $G(u)=G_{P}(u)\backslash \{\phi
(u)\}=X\backslash \{\phi (u)\}$.\medskip \newline
Note that $G\neq G_{P}$, as at every $u\in D^{\ast }$, $\phi (u)\in
G_{P}(u)\backslash G(u)$.\medskip

Observe that $G$ fails balancedness. Let $u$ be a profile in $D^{\ast }$
with $\phi (u)=x$, say. Let $y$ be the alternative ranked immediately above $%
x$ in \#1's ordering (nothing special about \#1 here).\ Some other
individual, say $i$, has $y$ ranked immediately below $x$. Construct $%
u^{\prime }$ from $u$ by transposition pair $(x,y)$ via 1 and $i$. Profile $%
u^{\prime }\in D$, so $G(u)=G_{P}(u)=x\neq G(u^{\prime })=X\backslash \{x\}$%
.\medskip

That this correspondence $G$ satisfies tops-in and Pareto is easy to see. \
Strong stability can be established by a detailed case-by-case
analysis.\medskip

Now that we have defined strong stability we can alter Theorem 3 by
substituting stability for monotonicity.\medskip

\textbf{Theorem 5}. \ For $m=4$ and $n\geq 2$, let $G:L(X)^{N}\rightarrow
2^{X}\backslash \{\varnothing \}$ be a social choice correspondence
satisfying all of: \medskip

\qquad 1. \ The Pareto condition;

\qquad 2. \ Tops-in;

\qquad 3. \ Balancedness;

\qquad 4. \ Strong stability;\medskip \newline
\relax then $G=G_{P}$, the Pareto correspondence. \medskip \newline
Proof: \ Suppose $w=\{a,b,c,w\}$ so that $m=4$ and that $G$ is a
correspondence satisfying the Pareto condition, tops-in, balancedness and
the strong stability condition. We need to show that if alternative $w$ is
Pareto-optimal at $u$, then it is in $G(u)$. \ Suppose otherwise, i.e.,
there exists a profile $u\in L(w)^{N}$ for which $w$ is Pareto optimal at $u$
but $w\notin G(u)$. \ Then $H(G)$ is defined. \ We show that each possible
value of $H(G)$ leads to a contradiction.\medskip

1. \ $H(G)=1$ is ruled out by tops-in.\medskip

2. \ $H(G)=4$ is ruled out by the Pareto condition.\medskip

3. \ Suppose $H(G)=3$. \ Then with a possible renumbering of individuals and
relabeling of alternatives, we have $n\geq 3$ and\medskip

\qquad \qquad $u$: \ 
\begin{tabular}{|l|l|l|l|}
\hline
$1$ & $2$ & $3$ & $\cdots $ \\ \hline
&  &  &  \\ 
$\vdots $ & $\vdots $ & $\vdots $ &  \\ 
$w$ & $w$ & $w$ & $\cdots $ \\ 
$a$ & $b$ & $c$ &  \\ \hline
\end{tabular}%
\medskip \newline
The alternative adjacent to but just above $w$ for \#1 is either $b$ or $c$.
\ But then a height argument shows a violation of $H(G)=3$.\medskip

4. \ Finally, suppose $H(G)=2$. \ There exists a profile $u$ with $w$ in say
\#1's second rank (with $a\neq w$ top-ranked), $w$ Pareto optimal at $u$ but 
$w\notin G(u).$\medskip

\qquad \qquad $u$: \ 
\begin{tabular}{|l|l|l|}
\hline
$1$ & $2$ & $\cdots $ \\ \hline
$a$ &  &  \\ 
$w$ &  &  \\ 
$\vdots $ &  & $\cdots $ \\ 
&  &  \\ \hline
\end{tabular}%
\medskip \newline
Because $w$ is Pareto optimal, someone, say \#2, must prefer $w$ to $a$. \
And $a$ can't be just below $w$ for \#2, for if it were, we could construct $%
u^{\ast }$ by transposing $a$ and $w$ for \#1 and \#2, resulting in $w$ at
\#1's top but, by balancedness, $w\notin G(u^{\ast })$, a violation of
tops-in. \ Without loss of generality, \#2 has $b$ between $w$ and $a$%
:\medskip

\qquad \qquad $u$: \ 
\begin{tabular}{|l|l|l|}
\hline
$1$ & $2$ & $\cdots $ \\ \hline
$a$ & $c$ &  \\ 
$w$ & $w$ &  \\ 
& $b$ & $\cdots $ \\ 
& $a$ &  \\ \hline
\end{tabular}%
\medskip \newline
Now $c$ can not be adjacent to but just below $w$ for \#1, or a
transposition of $w$ and $c$ for \#1 and \#2 would generate a profile with $%
w $ at \#2's top but $w$ not chosen. \ So we must have:\medskip

\qquad \qquad $u$: \ 
\begin{tabular}{|l|l|l|}
\hline
$1$ & $2$ & $\cdots $ \\ \hline
$a$ & $c$ &  \\ 
$w$ & $w$ &  \\ 
$b$ & $b$ & $\cdots $ \\ 
$c$ & $a$ &  \\ \hline
\end{tabular}%
\medskip \newline
But then look at profile\medskip

\qquad \qquad $u^{\prime }$: \ 
\begin{tabular}{|l|l|l|}
\hline
$1$ & $2$ & $\cdots $ \\ \hline
$a$ & $c$ &  \\ 
$w$ & $w$ &  \\ 
$c$ & $b$ & $\cdots $ \\ 
$b$ & $a$ &  \\ \hline
\end{tabular}%
\medskip \newline
Alternative $c\in G(u^{\prime })$ by tops-in. \ Profile $u$ can be obtained
from $u^{\prime }$ by lowering chosen $c$ just below $b$ for \#1. \ By
strong stability, this can not cause $w$ be kicked out of $G(u^{\prime })$.
\ But $w\notin G(u)$, so we must have had $w\notin G(u^{\prime })$. \ But
then a height argument shows a contradiction with $H(G)=2$. \ \ \ \ $\square 
$\medskip

\section{Final remarks\protect\medskip}

A perhaps unusual feature of the results in this paper is the use of the
Pareto condition in characterizing the Pareto social choice correspondence.
\ We make four observations regarding this:\medskip

1. \ What might at first appear to be an excessively strong condition is 
\textit{not} sufficient for five or more alternatives even when supplemented
with balancedness, tops-in, anonymity, and neutrality.\medskip

2. \ As noted earlier, the condition of excluding dominated alternatives is
actually \textit{extremely weak} and has been used in characterizing a wide
variety of standard social choice correspondences, e.g., Borda, plurality
voting, and the Copeland correspondence.\medskip

3. \ It helps to compare with an analogous use of a plurality condition. \
Suppose that we wanted to characterize plurality voting by using a \textit{%
plurality condition} that excludes all alternatives that are not plurality
winners and then added enough additional conditions to ensure that all
plurality winners \textit{are} included. \ This plurality condition would
seem quite artificially constructed, solely for the purpose of the one
characterization theorem. That's quite different from the Pareto condition.
\ Almost all standard social choice correspondences fail the plurality
condition.\footnote{%
One reader has suggested that we can simply characterize the Pareto
correspondence as the coarsest rule satisfying the Pareto condition. \ Of
course, in the same way, the plurality correspondence is the coarsest rule
satisfying the plurality condition; and the Borda correspondence is the
coarsest rule satisfying a Borda condition. \ None of these
"characterizations" really involve social choice properties (like
monotonicity, tops-in, stability, etc.) that relate social outcomes to
individual preferences. \ "Coarseness" is not a social choice property.}%
\medskip

4. \ Something like the Pareto condition \textit{is} required. \ Without
Pareto, we can get rules far from $G_{P}$: As seen by $G(u)=X$, \textit{all}
our other conditions combined, balancedness, monotonicity, and tops-in, plus
anonymity, neutrality, and strong stability are insufficient to rule out
correspondences that differ from $G_{P}$. \ And we can get rules very far
from $G(u)=X$ as well, even very close to $G_{P}$, but with just the
occasional choice of a dominated alternative, as seen next.\medskip

\textbf{Example 12}. \ Set $G(u)=G_{P}(u)$ except for those profiles with
complete agreement: $u(i)=u(j)$ (not just same top alternatives) for all $i$
and $j$. At profiles of complete agreement, set $G(u)$ to be the set
consisting of everyone's top \textit{two} alternatives (although the common
second is Pareto-dominated). \ This social choice correspondence satisfies
balancedness, monotonicity, and tops-in, plus anonymity, neutrality, and
strong stability.\medskip \bigskip

\textbf{REFERENCES} \medskip

Bellman, R (1984) \textit{Eye of the Hurricane}, Singapore: World
Scientific.\medskip

Campbell, DE, JS Kelly, and S Qi (2018): "A Stability Property in Social
Choice Theory"; International Journal of Economic Theory; Vol. 14.
85-95.\medskip

Kelly, JS and S Qi (2019): "Balancedness of Social Choice Correspondences."
arXiv: http://arxiv.org/abs/1804.02990, Mathematical Social Sciences,
forthcoming.\medskip

Moulin, H (1983) The Strategy of Social Choice (North-Holland).\medskip

Weymark, J (1984) "Arrow's Theorem with Social Quasi-Orderings"; Public
Choice; Vol. 42. 235-246.

\bigskip

\end{document}